\documentclass[11pt,lecno]{amsart}
\usepackage{indentfirst,amssymb,amsmath,amsthm}
\usepackage[mathscr]{euscript}
\usepackage[caption=false]{subfig}
\usepackage{graphicx}
\usepackage{tabulary}
\usepackage{epstopdf}

\newtheorem{theorem}{Theorem}[section]

\newtheorem{lemma}{Lemma}[section]

\newtheorem{corollary}{Corollary}[section]

\newcommand{\be}{\begin{equation}}
\newcommand{\ee}{\end{equation}}
\newcommand{\bea}{\begin{eqnarray}}
\newcommand{\eea}{\end{eqnarray}}
\newcommand{\beas}{\begin{eqnarray*}}
\newcommand{\eeas}{\end{eqnarray*}}

\begin{document}
\setcounter{page}{1} \setlength{\unitlength}{1mm}\baselineskip
.58cm \pagenumbering{arabic} \numberwithin{equation}{section}

\title[quasi-constant sectional curvature and $\mathcal{F}(\mathcal{R})$-gravity]
{Characterizations of a spacetime of quasi-constant sectional curvature and $\mathcal{F}(\mathcal{R})$-gravity }

\author[ U. C. De, K. De$^{*}$ , F. O. ZENGIN and  S. A. DEMIRBAG ]
{ Uday Chand De, Krishnendu De$^{*}$, Fusun OZEN ZENGIN and  Sezgin ALTAY DEMIRBAG}

\address
{Department of Pure Mathematics, University of Calcutta, West Bengal, India. ORCID iD: https://orcid.org/0000-0002-8990-4609}
\email {uc$_{-}$de@yahoo.com}
\address
 {$^{*}$ Department of Mathematics,
 Kabi Sukanta Mahavidyalaya,
The University of Burdwan.
Bhadreswar, P.O.-Angus, Hooghly,
Pin 712221, West Bengal, India. ORCID iD: https://orcid.org/0000-0001-6520-4520}
\email{krishnendu.de@outlook.in }

\address
{Department of Mathematics,  Istanbul Technical University, 34469,
Istanbul, Turkey. ORCID iD: https://orcid.org/ 0000-0002-5468-5100}
\email {fozen@itu.edu.tr}

\address
{Department of Mathematics,  Istanbul Technical University, 34469,
Istanbul, Turkey. ORCID iD: https://orcid.org/0000-0002-6643-6267}
\email {saltay@itu.edu.tr}

\footnotetext {$\bf{2020\ Mathematics\ Subject\ Classification\:}.$ 53B20; 53B30; 53B50; 83D05; 83C05.
\\ {Key words and phrases: Spacetime of quasi-constant sectional curvature; Perfect fluids; $f(\mathcal{R})$-gravity.\\
\thanks{$^{*}$ Corresponding author}
}}
\maketitle

\vspace{1cm}

\begin{abstract}

The main aim of this article is to investigate a spacetime of quasi-constant sectional curvature. At first, the existence of such a spacetime is established by several examples. We have shown that a spacetime of quasi-constant sectional curvature agrees with the present state of the universe and it represents a Robertson Walker spacetime. Moreover, if the spacetime is Ricci semi-symmetric or Ricci symmetric, then either the spacetime represents a spacetime of constant sectional curvature, or the spacetime represents phantom era. Also, we prove that a Ricci symmetric spacetime of quasi-constant sectional curvature represents a static spacetime and the spacetime under consideration is of Petrov type I, D or O. Finally, we concentrate on a quasi-constant sectional curvature spacetime solution in $\mathcal{F}(\mathcal{R})$-gravity. As a result, various energy conditions are studied and analysed our obtained outcomes in terms of a $\mathcal{F}(\mathcal{R})$-gravity model.

\end{abstract}

\maketitle

\section{Introduction}
To study a conformally flat hypersurfaces of a Euclidean space the authors \cite{chy} acquire the ensuing expression of the curvature tensor
\begin{eqnarray}\label{a1}
  R_{hijk}&=& \gamma(g_{hk}g_{ij}-g_{hj}g_{ik})\nonumber\\&&
  +\mu\{g_{hk}A_{i}A_{j}+g_{ij}A_{h}A_{k}-g_{hj}A_{i}A_{k}-g_{ik}A_{h}A_{j}\},
\end{eqnarray}
in which $\gamma, \mu$ are scalars and $A_{i}$ is a unit vector, called the generator.
An n-dimensional conformally flat space obeying (\ref{a1}) is named a space of quasi-constant sectional curvature and denoted by $(QC)_{n}$.\par 
However, if the equation (\ref{a1}) of the curvature tensor holds, then it can be easily verified that the space is conformally flat. So in the definition conformally flatness is not required. A space of quasi-constant sectional curvature has been studied by several authors, such as (\cite{boj}, \cite{chy} \cite{gan}, \cite{way}) and many others.\par

Spacetime is a time-oriented 4-dimensional Lorentzian manifold $\mathcal{M}$ which is a specific class of semi-Riemannian manifold endowed with a Lorentzian metric $g$ of signature $(-,+,+,+)$. A Lorentzian manifold is named a spacetime of quasi-constant sectional curvature if the curvature tensor $ R_{hijk}$ fulfills the condition (\ref{a1}). Here, we assume the generator $A_{i}$ is a unit time-like vector field, that is, $A_{i}A^{i}=-1$, $A^{i}=g^{ij}A_{j}$.\par

From the foregoing definition we can easily acquire
\begin{equation}\label{a2}
  R_{ij}=\{(n-1)\gamma+\mu\}g_{ij}+(n-2)\mu A_{i}A_{j};\;\;\;\;A_{i}A^{i}=-1,
\end{equation}
where $R_{hk}=R^{i}_{hki}$ is the Ricci tensor.\par
Multiplying the above equation with $g^{ij}$ yields
\begin{equation}\label{a3}
  \mathcal{R}=n(n-1)\gamma-2(n-1)\mu,
  \end{equation}
where $R=g^{hk}R_{hk}$ is the Ricci scalar.
Again multiplying (\ref{a2}) by $A^{i}A^{j}$, we get
\begin{equation}\label{a4}
  R_{ij}A^{i}A^{j}=-(n-1)\gamma+(n-3)\mu.
\end{equation}
Let us denote the Ricci curvature $R_{ij}A^{i}A^{j}$ by $\nu$.\par
Hence from equations (\ref{a3}) and (\ref{a4}), we obtain
\begin{equation}\label{a5}
  \gamma=\frac{(n-3)\mathcal{R}-2\nu}{(n-1)\{n(n-3)+2\}};\;\;\;\mu=\frac{\mathcal{R}+n\nu}{n(n-3)+2}.
\end{equation}
Therefore, the associated scalars $\gamma$ and $\mu$ are related with the Ricci scalar $\mathcal{R}$ and Ricci curvature $\nu$.\par

In a spacetime $\mathcal{M}^{n}$, for a smooth function $\psi>0$ (also called scale factor, or warping function), if $\mathcal{M}=-I \times_\psi,\mathrm{M}$, where $I$ is the open interval of $\mathbb{R}$ 
, $\mathrm{M}^{n-1}$ denotes the Riemannian manifold, then $\mathcal{M}$ is named a generalized Robertson Walker (briefly, GRW) spacetime \cite{alias1}. If $\mathcal {M}^{3}$ is of constant sectional curvature, then the spacetime represents a Robertson Walker (briefly, RW) spacetime.\par

Throughout this paper we consider 4-dimensional spacetime of quasi-constant sectional curvature, denoted by $(QC)_4$- spacetime.\par

Due to the absence of a stress tensor and heat conduction terms corresponding to viscosity, the fluid is referred to as perfect and the energy momentum tensor $T_{hk}$ is written by
\begin{equation}
\label{a9}
T_{hk}=(\sigma+p)u_{h}u_{k}+p g_{hk},
\end{equation}
in which $p$ and $\sigma$ stand for the perfect fluid's isotropic pressure and energy density, respectively \cite{o'neill}.  In the last equation, $g$ is the Lorentzian metric and $u_{h}$, the velocity vector is defined by $g_{hk}u^{h}u^{k}=-1$ and $u_{h}=g_{hk}u^{k}$.\par

For a gravitational constant $\kappa$, the Einstein's field equations (briefly, EFEs) without a cosmological constant is described by
\begin{equation}
\label{a10}
R_{hk}-\frac{R}{2}g_{hk}=\kappa T_{hk},
\end{equation}
where $R_{hk}=R^{i}_{hki}$ and $R=g^{hk}R_{hk}$ indicate the Ricci tensor and the Ricci scalar, respectively.\par

A spacetime $\mathcal{M}$ is named a perfect fluid (briefly PF) spacetime if the non-vanishing Ricci tensor $R_{hk}$ fulfills
\begin{equation}
\label{a11}
R_{hk}=\alpha g_{hk}+\beta u_{h} u_{k},
\end{equation}
in which $\alpha$ and $\beta$ are smooth functions. The foregoing equation is obtained from the equations (\ref{a9}) and (\ref{a10}).\par

Combining the equations (\ref{a9}), (\ref{a10}) and (\ref{a11}), we acquire
\begin{equation}
\label{a12}
\beta=k^2 (p+\sigma), \,\, \alpha=\frac{k^2 (p-\sigma)}{2-n}.
\end{equation}


In \cite{o'neill}, O'Neill established that a RW- spacetime is a PF- spacetime. Every GRW- spacetime of dimension four is a PF- spacetime if it is a RW- spacetime and vice-versa \cite{gtt}. For additional insights about the PF-spacetimes and GRW spacetimes, we refer( \cite{blaga2}, \cite{bychen}, \cite{de}, \cite{survey}) and its references.\par

Additionally, a state equation with the form $p = p(\sigma)$ connects $p$ and $\sigma$, and the PF-spacetime is known as isentropic. Furthermore, if $p = \sigma$, the PF- spacetime is referred to as stiff matter. The PF- spacetime is referred to as the dark matter era if $p+\sigma=0 $, the dust matter fluid if $p = 0$, and the radiation era if $p =\frac{\sigma}{3}$ \cite{ch1}. The universe is represented as accelerating phase when $\frac{p}{\sigma}< {-\frac{1}{3}}$. It covers the quintessence phase if $-1< \frac{p}{\sigma}< 0$ and phantom era if $\frac{p}{\sigma}< -1$.\par
In this paper, we consider a $(QC)_4$- spacetime and prove the following:
\begin{theorem}\label{th1}
A spacetime of quasi-constant sectional curvature represents a perfect fluid spacetime and agrees with the present state of the universe.
\end{theorem}
As a consequence, we also show that the $(QC)_4$- spacetime represents a RW- spacetime.\par

The Ricci tensor is said to be of Codazzi type \cite{gray} if it fulfills the subsequent relation
\begin{equation}\nonumber
 \nabla_{l} R_{hk}= \nabla_{k}R_{hl}.
\end{equation}
Here we discuss a $(QC)_{4}$-spacetime fulfilling the Ricci tensor is of Codazzi type and establish the subsequent theorem:
\begin{theorem}\label{th2}
A $(QC)_{4}$-spacetime obeying $\gamma^{h}A_{h}=0$ endowed with Codazzi type of Ricci tensor represents a spacetime of constant sectional curvature, or the vorticity vanishes.
\end{theorem}

A spacetime is named semi-symmetric \cite{sz} if it obeys the relation
\begin{equation}\label{a6}
 \nabla_{l}\nabla_{m} R^{h}_{ijk}-\nabla_{m}\nabla_{l}R^{h}_{ijk}=0,
\end{equation}
in which $\nabla$ indicates the covariant differentiation. It is to be noted that the class of locally symmetric spaces ($\nabla_{l} R^{h}_{ijk}=0$) due to Cartan is a proper subset of semi-symmetric spaces.\par

A spacetime is named Ricci semi-symmetric \cite{mi} if it fulfills the relation
\begin{equation}\label{a7}
  \nabla_{l}\nabla_{m}R_{ij}-\nabla_{m}\nabla_{l}R_{ij}=0.
\end{equation}
In this article we consider a Ricci semi-symmetric $(QC)_{4}$-spacetime and establish the following:
\begin{theorem}\label{th3}
If a $(QC)_{4}$-spacetime is Ricci semi-symmetric, then either the spacetime represents a spacetime of constant sectional curvature, or the spacetime represents phantom era.
\end{theorem}
Again, the class of Ricci symmetric spaces ($\nabla_{l} R_{ij}=0$) is a proper subset of Ricci semi-symmetric spaces. Every semi-symmetric space is known to be Ricci semi-symmetric, but the converse is usually not true. In a Riemannian space they are equivalent for dimension three. In \cite{ta}, it has been established that for $n\geq 3$, the above stated relations are equivalent for hypersurfaces having non negative scalar curvature in a Euclidean space $E^{n+1}$.\par
Inspired by the above studies we prove the subsequent result:
\begin{corollary}\label{th4}
A Ricci symmetric $(QC)_{4}$-spacetime either represents a spacetime of constant sectional curvature, or the phantom era.
\end{corollary}

If a Lorentzian manifold permits a timelike Killing vector field $\rho$, it is referred to as a stationary spacetime and static (\cite{sanches}, \cite{ste}, p. 283) if, additionally, $\rho$ is irrotational. We will refer to $\rho$ in this context as the static vector field, where it is assumed that spacetime is time-oriented. The product $\mathbb{R}\times S$ is called a static spacetime if it is equipped with the metric
\begin{equation}\label{sta}
  g[(t, y)]=-\beta(y)dt^2 + g_{S}[y],
\end{equation}
in which $g_{S}$ denotes a Riemannian metric on S. Any static spacetime behaves like a standard one locally, with $\rho$ recognisable to $\partial t$.\par

Here we also consider a Ricci symmetric $(QC)_{4}$-spacetime and prove the following:
\begin{theorem}\label{th5}
If a $(QC)_{4}$-spacetime with $\mu \neq 0$ is Ricci symmetric, then the spacetime represents a static spacetime and the spacetime is of Petrov type I, D or O.
\end{theorem}

EFEs are incapable to answer for the universe's late-time inflation without accepting the existence of dark energy. This inspired other scholars to expand it in order to obtain higher order gravity field equations. One of the aforementioned modified theories of gravity was developed by substituting an arbitrary function $\mathcal{F}(\mathcal{R})$ for the Ricci scalar $\mathcal{R}$ in the Einstein-Hilbert action, which was first suggested in 1970 by H.A. Buchdahl \cite{hab}. Of course, various theoretical scalar-tensor findings and observational data, constrain the feasibility of such functions. Many useful functional versions of $\mathcal{F}(\mathcal{R})$ theory have been developed recently; for additional information, see (\cite{cap1}, \cite{harfr}, \cite{ade1}, \cite{kde}, \cite{kde1}). By the motivation of these investigations, in this article we concentrate on a $(QC)_4$- spacetime solution in $\mathcal{F}(\mathcal{R})$-gravity and prove the following main results.\par
\begin{theorem}\label{th6}
For $\mathcal{R}=$ constant, in any $(QC)_4$ spacetime solution of the $\mathcal{F}(\mathcal{R})$-gravity, the Ricci tensor $R_{ij}$ is semi-symmetric if and only if the energy momentum tensor $T_{ij}$ is semi-symmetric.
\end{theorem}
\begin{theorem}\label{th7}
For $\mathcal{R}=$ constant, in any $(QC)_4$ spacetime solution of the $\mathcal{F}(\mathcal{R})$-gravity, if the Ricci tensor $R_{ij}$ is symmetric then either the spacetime represents dark matter, or the energy density is constant.
\end{theorem}
\begin{theorem}
\label{th8}
For $\mathcal{R}=$ constant, in any $(QC)_4$ spacetime solution of the $\mathcal{F}(\mathcal{R})$-gravity, the isotropic pressure $p$ and energy density $\sigma$ obey the subsequent equations
$$p=\frac{(3\gamma-\mu)\mathcal{F}_\mathcal{R}(\mathcal{R})}{\kappa^{2}}
-\frac{\mathcal{F}(\mathcal{R})}{2\kappa^{2}},$$
$$\sigma=\frac{3(\mu-\gamma)\mathcal{F}_\mathcal{R}(\mathcal{R})}{\kappa^{2}}
+\frac{\mathcal{F}(\mathcal{R})}{2\kappa^{2}}.$$
\end{theorem}
The energy conditions (ECs), as is well known in GR, are essential tools for studying wormholes and black holes in various modified gravities (\cite{bam}, \cite{hak}). The ECs are carefully constructed using the Raychaudhuri equations \cite{ray}, which reflect the attractiveness of gravity through the positive condition $\mathcal{R}_{lk}v^{l}v^{k}\geq0,$ where $v^{l}$ is a null vector. This condition on geometry in GR is comparable to the null energy condition (NEC) $\mathcal{T}_{lk}v^{l}v^{k}\geq0$ on matter through the EFE. In particular, the weak energy condition (WEC), states that $\mathcal{T}_{lk}u^{l}u^{k}\geq0,$ for every timelike vector $u^{l}$ and under the assumption of a positive local energy density.\par
ECs are coordinate-invariant constraints on the energy-momentum tensor (effective) that are helpful when considering the probability of different matter origins, not just a PF continuum, and which fulfill the theories of modified gravity while maintaining the notion that energy has to be positive. There are various ECs; some of them, like the trace EC, are no longer relevant today, while others are weaker and combined with others. However, the main idea is to construct various scalar fields by contracting the energy momentum tensor with any lightlike ,or timelike  vectors. For our topic, we derive several important ECs using our obtained results. Finally, we demonstrate some ECs by considering a $\mathcal{F}(\mathcal{R})$ model.

\section{Examples of $(QC)_4$- spacetime}
\subsection{Example 1\;:}
In \cite{survey}, Mantica and Molinari established that every RW-spacetime is a $(QC)_4$- spacetime.
\subsection{Example 2\;:}
The pseudo-symmetric manifold established by Chaki \cite{cha}, whose Riemann curvature fulfills the condition
\begin{equation*}
  \nabla_{l}R_{hijk}=2A_{l}R_{hijk}+A_{h}R_{lijk}+A_{i}R_{hljk}+A_{j}R_{hilk}+A_{k}R_{hijl},
\end{equation*}
in which $A_{i}$ indicates the non-zero 1-form. In \cite{tarafdar1}, Tarafdar established that a conformally flat pseudo-symmetric manifold is a space of quasi constant curvature.
\subsection{Example 3\;:}
In a 4-dimensional conformally flat spacetime, $R_{hijk}$ (curvature tensor) is given by
 \begin{eqnarray}\label{ex1}
  R_{hijk}&=& \frac{1}{2}(g_{hk}R_{ij}-g_{hj}R_{ik}+g_{ij}R_{hk}-g_{ik}R_{hj})\nonumber\\&&
  +\frac{\mathcal{R}}{6}\{g_{hk}g_{ij}-g_{hj}g_{ik}\},
\end{eqnarray}
A semi-Riemannian manifold $M$ is said to be Ricci simple \cite{dep} if $$R_{ij}=-\mathcal{R}A_{i}A_{j},$$
$A_{i}$ is a unit time-like vector. The physical interpretation of Ricci simple spacetime is given by Mantica and Molinari \cite{survey} in which they prove that it represents stiff matter fluid. Now using the above relation in \ref{ex1} we obtain that the curvature tensor is of the form of a $(QC)_4$. Hence, every Ricci simple conformally flat spacetime is a $(QC)_4$- spacetime.

\subsection{Example 4\;:}
A weakly Ricci symmetric manifold invented by Tamassy and Binh \cite{tam}, whose Ricci tensor $R_{ij}$ is
not identically zero and obeys the condition
\begin{equation*}
  \nabla_{l}R_{ij}=A_{l}R_{ij}+B_{i}R_{lj}+C_{j}R_{il},
\end{equation*}
in which $A, B, C$ indicate the non-zero 1-forms.\par
In a $(WRS)_{4}$ (weakly Ricci symmetric spacetime) \cite{desk}, $R_{ij}=-\mathcal{R}u_{i}u_{j}.$ Using this in \ref{ex1} we arrive at the same result. Therefore, every conformally flat $(WRS)_{4}$ spacetime is a $(QC)_4$- spacetime.

\section{Spacetime of quasi-constant sectional curvature }
{\bf Proof of the Theorem \ref{th1} :}\par
Multiplying (\ref{a1}) with $g^{ij}$, we acquire
\begin{eqnarray}\label{b1}
  R_{hk} &=& \gamma(4g_{hk}-g_{hk})+\mu\{-g_{hk}+4A_{h}A_{k}-A_{h}A_{k}-A_{h}A_{k}\} \nonumber\\&&
  =(3\gamma-\mu)g_{hk}+2\mu A_{h}A_{k},
\end{eqnarray}
which represents a PF- spacetime.\par
Multiplying (\ref{b1}) with $g^{hk}$, we obtain
\begin{equation}\label{b2}
  \mathcal{R}=6(2\gamma-\mu).
\end{equation}
Comparing the previous equation with (\ref{a12}), we get

\begin{equation}
\label{b3}
2\mu=k^2 (p+\sigma),
\end{equation}
and
\begin{equation}
\label{b4}
 3\gamma-\mu=-\frac{k^2 (p-\sigma)}{2}.
\end{equation}
From the last two equations, we infer
\begin{equation}\label{b5}
  p=\frac{1}{\kappa^{2}}(-3\gamma+2\mu)\;\;\; \sigma=\frac{1}{\kappa^{2}}3\gamma.
\end{equation}
Since $p$ and $\sigma$ are not constant, we can say the $(QC)_4$- spacetime agrees with the present state of the universe.\par
Hence, the theorem is proved.\par
\vspace{.2cm}
If $3\gamma=2\mu$, the equation (\ref{b5}) gives $p=0$. Hence, we state
\begin{corollary}
  A spacetime of quasi-constant sectional curvature with $3\gamma=2\mu$ represents dust era.
\end{corollary}
From equation (\ref{b5}), we acquire $p+\sigma=\frac{2\mu}{\kappa^{2}}\neq 0$. Also, we know that the space of quasi-constant sectional curvature is conformally flat and hence $div C=0$ (`$div$' denotes the divergence), where C is the Weyl curvature tensor. In \cite{survey}, it is established that a PF- spacetime with $p+\sigma \neq 0$  and $div\, C = 0$, is a GRW- spacetime.\par

Again, for $n=4$, every GRW- spacetime is a PF- spacetime if and only if it is a RW- spacetime \cite{gtt} and thus the spacetime becomes RW- spacetime.
\begin{corollary}
\label{rem2.2}
A spacetime of quasi-constant sectional curvature represents a RW- spacetime.
\end{corollary}

{\bf Proof of the Theorem \ref{th2} :}\par
Let the Ricci tensor be of Codazzi type, which entails
\begin{equation}\label{1b}
 \nabla_{l} R_{hk}= \nabla_{k}R_{hl}.
\end{equation}
Differentiating (\ref{b1}) covariantly yields
\begin{equation}\label{c1}
  \nabla_{l}R_{hk}=(3\gamma_{l}-\mu_{l})g_{hk}+2\mu_{l}A_{h}A_{k}+2\mu(\nabla_{l}A_{h}A_{k}-A_{h}\nabla_{l}A_{k}).
\end{equation}
Similarly, we get
\begin{equation}\label{c2}
  \nabla_{k}R_{hl}=(3\gamma_{k}-\mu_{k})g_{hl}+2\mu_{k}A_{h}A_{l}+2\mu(\nabla_{k}A_{h}A_{l}-A_{h}\nabla_{k}A_{l}).
\end{equation}
 Hence using the foregoing equations in (\ref{1b}), we acquire
\begin{eqnarray}\label{c3}
  0&=&\nabla_{l}R_{hk}-\nabla_{k}R_{hl}\nonumber\\&&
   = (3\gamma_{l}-\mu_{l})g_{hk}-(3\gamma_{k}-\mu_{k})g_{hl}\nonumber\\&&
  +2\mu_{l}A_{h}A_{k}-2\mu_{k}A_{h}A_{l}+2\mu(\nabla_{l}A_{h}A_{k}-A_{h}\nabla_{l}A_{k})
  -2\mu(\nabla_{k}A_{h}A_{l}-A_{h}\nabla_{k}A_{l}).
\end{eqnarray}
Multiplying the above equation by $g^{hk}$ gives
\begin{eqnarray}\nonumber
  0&=& 4(3\gamma_{l}-\mu_{l})-(3\gamma_{l}-\mu_{l})\nonumber\\&&
  -2\mu_{l}-2\mu^{h}A_{h}A_{l}-2\mu(\nabla_{h}A^{h}A_{l}+A^{k}\nabla_{k}A_{l})\nonumber,
\end{eqnarray}
which implies
\begin{eqnarray}\label{c4}
  0&=& 9\gamma_{l}-5\mu_{l}-2\mu^{h}A_{h}A_{l}\nonumber\\&&
  -2\mu(\nabla_{h}A^{h}A_{l}+A^{k}\nabla_{k}A_{l}).
\end{eqnarray}
Equation (\ref{1b}) implies $\mathcal{R}=$ constant and hence from (\ref{b2}), we acquire $2\gamma_{l}=\mu_{l}$. Using this in the above equation we infer
\begin{equation}\label{c5}
  0= -\gamma_{l}-4\gamma^{h}A_{h}A_{l} -2\mu(A^{h}_{h}A_{l}+A^{k}\nabla_{k}A_{l}).
\end{equation}

Multiplying (\ref{c5}) by $A^{l}$, we obtain
\begin{equation}\nonumber
  0= -\gamma_{l}A^{l}+4\gamma^{h}A_{h} -2\mu A^{h}_{h},
\end{equation}
which implies
\begin{equation}\label{c5}
  0= 3\gamma^{h}A_{h} -2\mu A^{h}_{h}.
\end{equation}
Let us consider $\gamma^{h}A_{h}=0$. Then the previous equation entails that either $\mu=0$, or $div A^{h}=0$.\par
If $\mu=0$, the spacetime of quasi constant curvature reduces to the spacetime of constant sectional curvature.\par
If $div A^{h}=0$, then the velocity vector field is conservative. A conservative vector field is always irrotational, thus we conclude that the PF has zero vorticity.\par
This ends the proof.\par

{\bf Proof of the Theorem \ref{th3} :}\par
Let a $(QC)_{4}$-spacetime be Ricci semi-symmetric, that is
\begin{equation}
\label{aa7}
 \nabla_{l}\nabla_{m} R_{ij}-\nabla_{m}\nabla_{l}R_{ij}=0.
\end{equation}
Hence using (\ref{b1}) in the above equation, we acquire
\begin{equation}\label{bb2}
  \mu\{\nabla_{l}\nabla_{m}A_{h}A_{k}+A_{h}\nabla_{l}\nabla_{m}A_{k}
  -\nabla_{m}\nabla_{l}A_{h}A_{k}-A_{h}\nabla_{m}\nabla_{l}A_{k}\}=0,
\end{equation}
This implies either $\mu=0$, or $\mu\neq 0$.\par 

Case (i): If $\mu=0$, the spacetime of quasi constant curvature reduces to the spacetime of constant curvature.\par

Case (ii): If $\mu\neq 0$, then $(\nabla_{l}\nabla_{m}A_{h}-\nabla_{m}\nabla_{l}A_{h})A_{k}+A_{h}(\nabla_{l}\nabla_{m}A_{k}
  -\nabla_{m}\nabla_{l}A_{k})=0$. Hence applying Ricci identity, we get
\begin{equation}\label{bb3}
  \{A_{p}R^{p}_{hlm}A_{k}+A_{h}A_{p}R^{p}_{klm}\}=0.
\end{equation}

Multiplying with $A^{k}$ the foregoing equation yields
\begin{equation}\label{bb4}
 \{-A_{p}R^{p}_{hlm}+A^{p}A^{k}R_{pklm}\}=0.
\end{equation}
Since $R_{pkml}=-R_{kpml}$, therefore $A^{p}A^{k}R_{pklm}=0.$\par

Hence the above equation gives
\begin{equation}\label{bb5}
 A_{p}R^{p}_{hlm}=0,
\end{equation}
which implies
\begin{equation}\label{bb6}
 A_{p} R^{p}_{h}=0.
\end{equation}
Now from (\ref{b1}), we obtain
\begin{equation}\label{bb7}
  A^{k}R_{hk} =(3\gamma-\mu)A_{h}-2\mu A_{h}.
\end{equation}
Using (\ref{bb6}) we can write
\begin{equation}\label{bb8}
  (3\gamma-\mu)A_{h}-2\mu A_{h}=0,
\end{equation}
which implies $\mu=\gamma$.\par
Finally, we get either $\mu=0$, or $\mu=\gamma$.\par
If $\mu=\gamma$, then using this relation in (\ref{b5}), we infer $\frac{p}{\sigma}=-\frac{1}{3}$ which means the spacetime represents phantom era.\par
This finishes the proof.

{\bf Proof of the Corollary \ref{th4} :}\par
Since Ricci symmetry ($\nabla_{l}R_{ij}=0$) implies Ricci semi-symmetry ($\nabla_{l}\nabla_{m} R_{ij}-\nabla_{m}\nabla_{l}R_{ij}=0$), therefore from the above result for a Ricci symmetric $(QC)_{4}$-spacetime we can state either the spacetime represents a spacetime of constant sectional curvature, or the spacetime represents phantom era.\par
Therefore the proof is completed.

{\bf Proof of the Theorem \ref{th5} :}\par
Let  $\mu\neq 0$ and $\mu=\gamma$.\par
Using $\mu=\gamma$ in (\ref{b1}), we acquire
\begin{equation}\label{r1}
  R_{ij}=2\mu(g_{ij}+A_{i}A_{j}).
\end{equation}
Covariant differentiation of the foregoing equation yields
\begin{equation}\label{r2}
  \nabla_{l}R_{ij}=2\nabla_{l}\mu(g_{ij}+A_{i}A_{j})+2\mu(\nabla_{l}A_{i}A_{j}+A_{i}\nabla_{l}A_{j}).
\end{equation}
By hypothesis $\nabla_{l}R_{ij}=0$. Hence from the above we infer
\begin{equation}\label{r3}
\nabla_{l}\mu(g_{ij}+A_{i}A_{j})+\mu(\nabla_{l}A_{i}A_{j}+A_{i}\nabla_{l}A_{j})=0.
\end{equation}
Multiplying by $A^{i}$ gives
\begin{equation}\label{r4}
  \mu \nabla_{l}A_{j}=0,
\end{equation}
since $A^{i}\nabla_{l}A_{i}=0$.\par

Therefore, we obtain $\nabla_{l}A_{j}=0$, since $\mu \neq 0$.\par

Now $\nabla_{l}A_{j}=0$ implies that $A_{j}$ is a Killing vector and also irrotational. Hence the spacetime becomes static \cite{ste}.\par

It is well-known (\cite{ehl}, Section 10.7) that any static spacetime is everywhere of Petrov type I, D or O. As a result, the spacetime under consideration is of Petrov type I, D or O.\par
This ends the proof.

\section{ $\mathcal{F}(\mathcal{R})$-Gravity}


Taking the Einstein-Hilbert action term into account
\begin{align}
S=\frac{1}{2\kappa^2}\int{{\sqrt{-g}}\mathcal{F}(\mathcal{R})d^4x}+\int{{\sqrt{-g}}L_{m}d^4x}
\nonumber
\end{align}
where 
$k=\sqrt{8{\pi}G}$, $G$ is Newton's constant and $L_m$ is referred as the matter Lagrangian density and described by
\begin{align}
T_{ij}=-\frac{2}{\sqrt{-g}}\frac{\delta (\sqrt{-g}L_m)}{\delta g^{ij}}
\nonumber
\end{align}
The field equation by utilizing the variation with regard to 
$g^{ij}$, can have the form

\begin{align}\label{g1}
\mathcal{F}_\mathcal{R}(\mathcal{R})R_{ij}-(\nabla_i\nabla_j\mathcal{F}_\mathcal{R}(\mathcal{R})-g_{ij}{\square})
-\frac{\mathcal{F}(\mathcal{R})}{2}g_{ij}=k^2T_{ij}
\end{align}
where the D'Alembert operator is represented by $\square$ and the differentiation with respect to $\mathcal{R}$ is denoted by "$\mathcal{F}_\mathcal{R}(\mathcal{R})$".
If $\frac{\mathcal{R}F_\mathcal{R}(\mathcal{R})}{2} g_{ij}$ is subtracted from both sides of \eqref{g1}, then the result is
\begin{align}\label{g2}
\mathcal{F}_\mathcal{R}(\mathcal{R})R_{ij}-\frac{\mathcal{R}F_\mathcal{R}(\mathcal{R})}{2}g_{ij}
=k^2T_{ij}+k^2T^{(curve)}_{ij}
\end{align}

such that
\begin{align}\label{g3}
T^{(eff)}_{ij}=T_{ij}+T^{(curve)}_{ij}
\end{align}
where
\begin{align}\label{g4}
T^{(curve)}_{ij}=\frac{1}{k^2}\big[(\nabla_i\nabla_j-g_{ij}{\square})\mathcal{F}_\mathcal{R}(\mathcal{R})+\frac{(\mathcal{F}(\mathcal{R})
-\mathcal{R}\mathcal{F}_\mathcal{R}(\mathcal{R}))}{2}g_{ij}
\big].
\end{align}
$T^{(eff)}_{ij}$ represents an efficient voltage-energy tensor.

With the help of \eqref{g2} and \eqref{g4} in case of constant $\mathcal{R}$, we acquire
\begin{align}\label{g5}
R_{ij}-\frac{\mathcal{R}}{2}g_{ij}=\frac{\kappa^2}{\mathcal{F}_\mathcal{R}(\mathcal{R})}T_{ij}
+\big[\frac{(\mathcal{F}(\mathcal{R})-\mathcal{R}\mathcal{F}_\mathcal{R}(\mathcal{R})}{2\mathcal{F}_\mathcal{R}(\mathcal{R})}g_{ij}\big].
\end{align}	

\begin{lemma}
\label{l1}
\cite{ade2} For $\mathcal{R}=$ constant, in a PF- spacetime solution of the $\mathcal{F}(\mathcal{R})$-gravity, the effective pressure $p^{(eff)}$ and the effective energy density $\sigma^{(eff)}$ obey the subsequent equations
$$p^{(eff)}=p+\frac{1}{2\kappa^2}\big[(\mathcal{F}(\mathcal{R})
-\mathcal{R}\mathcal{F}_\mathcal{R}(\mathcal{R})\big],$$
$$\sigma^{(eff)}=\sigma-\frac{1}{2\kappa^2}\big[\mathcal{F}(\mathcal{R})
-\mathcal{R}\mathcal{F}_\mathcal{R}(\mathcal{R})\big].$$
\end{lemma}

In Theorem \ref{th1}, we establish that every $(QC)_4$ spacetime is a PF-spacetime. Hence, the energy momentum tensor is of the form (\ref{a9}). Here we assume that the velocity vector field $u_{i}$ of the PF-spacetime is identical with the associated vector $A_{i}$ of the $(QC)_4$ spacetime. \par
Now, using (\ref{a9}) in (\ref{g5}), we acquire
\begin{align}\label{h1}
R_{ij}-\frac{\mathcal{R}}{2}g_{ij}=\frac{\kappa^2}{\mathcal{F}_\mathcal{R}(\mathcal{R})}\{(p+\sigma)A_{i}A_{j}+pg_{ij}\}
+\big[\frac{(\mathcal{F}(\mathcal{R})-\mathcal{R}\mathcal{F}_\mathcal{R}(\mathcal{R})}{2\mathcal{F}_\mathcal{R}(\mathcal{R})}g_{ij}\big].
\end{align}	
Multiplying both sides of the previous equation by $g^{ij}$ gives
\begin{align}\label{h2}
\mathcal{R}-2\mathcal{R}=\frac{\kappa^2}{\mathcal{F}_\mathcal{R}(\mathcal{R})}\{-(p+\sigma)+4p\}
+2\big[\frac{(\mathcal{F}(\mathcal{R})-\mathcal{R}\mathcal{F}_\mathcal{R}(\mathcal{R})}{\mathcal{F}_\mathcal{R}(\mathcal{R})}\big],
\end{align}	
which implies
\begin{align}\label{h3}
\mathcal{R}=\frac{1}{\mathcal{F}_\mathcal{R}(\mathcal{R})}\{\kappa^{2}(3p-\sigma)+2\mathcal{F}(\mathcal{R})\}.
\end{align}
Putting the value of $\mathcal{R}$ in (\ref{h1}) yields
\begin{align}\label{h4}
R_{ij}=\frac{2\kappa^2 p+\mathcal{F}(\mathcal{R})}{2\mathcal{F}_\mathcal{R}(\mathcal{R})} g_{ij}+\frac{\kappa^2 (p+\sigma)}{\mathcal{F}_\mathcal{R}(\mathcal{R})}A_{i}A_{j}.
\end{align}
\begin{theorem}
For $\mathcal{R}=$ constant, in any $(QC)_4$ spacetime solution of the $\mathcal{F}(\mathcal{R})$-gravity theory, the Ricci tensor $R_{ij}$ is of the form (\ref{h4}).
\end{theorem}


{\bf Proof of the Theorem \ref{th6} :}\par
From the field equation (\ref{h1}), we can easily acquire
\begin{align}\label{h6}
\nabla_{l}\nabla_{m}T_{ij}-\nabla_{m}\nabla_{l}T_{ij}=0 \Leftrightarrow \nabla_{l}\nabla_{m}R_{ij}-\nabla_{m}\nabla_{l}R_{ij}=0.
\end{align}
Hence, the theorem is proved.\par

{\bf Proof of the Theorem \ref{th7} :}\par
Covariant differentiation of the equation (\ref{h4}) yields
\begin{align}\label{h7}
\nabla_{l}\nabla_{m}R_{ij}= \nabla_{l}\nabla_{m}a g_{ij}+\nabla_{l}\nabla_{m} b A_{i}A_{j}
+b\{\nabla_{l}\nabla_{m} A_{i}A_{j}+A_{i}\nabla_{l}\nabla_{m}A_{j}\}.
\end{align}
where $a=\frac{2\kappa^2 p+\mathcal{F}(\mathcal{R})}{2\mathcal{F}_\mathcal{R}(\mathcal{R})}$ and $b=\frac{\kappa^2 (p+\sigma)}{\mathcal{F}_\mathcal{R}(\mathcal{R})}$.
The foregoing equation immediately gives
\begin{align}\label{h8}
\nabla_{l}\nabla_{m}R_{ij}-\nabla_{m}\nabla_{l}R_{ij}=b\{ \nabla_{l}\nabla_{m} A_{i}A_{j}-\nabla_{m}\nabla_{l}A_{i}A_{j}+A_{i}\nabla_{l}\nabla_{m}A_{j}-A_{i}\nabla_{m}\nabla_{l}A_{j}\}.
\end{align}
Let the Ricci tensor be semi-symmetric, that is, $R^{h}_{ijk,lm}-R^{h}_{ijk,ml}=0$. Then using Ricci identity the above equation infers that
\begin{align}\label{h9}
0=b\{ A_{j}A_{h}R^{h}_{ilm}+A_{i}A_{h}R^{h}_{jlm}\},
\end{align}
which implies either $b=0$, or $b\neq0.$\par
Case (i): If $b=0$, then we have $p+\sigma=0$. Therefore, the spacetime represents the dark matter era.\par
Case (ii): If $b\neq0$, then $A_{j}A_{h}R^{h}_{ilm}+A_{i}A_{h}R^{h}_{jlm}=0.$ Therefore multiplying the previous relation by $g^{lm}$ we have
\begin{align}\label{h10}
A_{j}A_{h}R^{h}_{i}+A_{i}A_{h}R^{h}_{j}=0.
\end{align}
Also, from (\ref{h4}), we acquire
\begin{align}\label{h11}
A^{j}R_{ij}=aA_{i}-bA_{i}=(a-b)A_{i}.
\end{align}
Using (\ref{h11}) in (\ref{h10}), we get
\begin{align}\label{h12}
(a-b)A_{i}A_{j}=0.
\end{align}
From the previous equation we conclude that $a=b$ which implies $\sigma=\frac{\mathcal{F}(\mathcal{R})}{2\kappa^2}=$ constant.\par
This ends the proof.\par

{\bf Proof of the Theorem \ref{th8} :}\par
Using (\ref{b1}) the field equation turns into
\begin{align}\label{g6}
(3\gamma-\mu)g_{ij}+2\mu A_{i}A_{j}-\frac{\mathcal{R}}{2}g_{ij}=\frac{\kappa^2}{\mathcal{F}_\mathcal{R}(\mathcal{R})}T_{ij}
+\big[\frac{(\mathcal{F}(\mathcal{R})-\mathcal{R}\mathcal{F}_\mathcal{R}(\mathcal{R})}{2\mathcal{F}_\mathcal{R}(\mathcal{R})}g_{ij}\big].
\end{align}	

Multiplying both sides of the foregoing equation by $g^{ij}$ yields
\begin{align}\label{g7}
4(3\gamma-\mu)-2\mu =\frac{\kappa^2}{\mathcal{F}_\mathcal{R}(\mathcal{R})}\big[3p-\sigma\big]+\frac{2\mathcal{F}(\mathcal{R})}{\mathcal{F}_\mathcal{R}(\mathcal{R})}.
\end{align}	
Again multiplying (\ref{g6}) by $A^{i}A^{j}$, we obtain
\begin{align}\label{g8}
-(3\gamma-\mu)+2\mu =\frac{\kappa^2 \sigma}{\mathcal{F}_\mathcal{R}(\mathcal{R})}-\frac{\mathcal{F}(\mathcal{R})}{2\mathcal{F}_\mathcal{R}(\mathcal{R})},
\end{align}	
which implies
\begin{align}\label{g9}
\sigma=\frac{3(\mu-\gamma)\mathcal{F}_\mathcal{R}(\mathcal{R})}{\kappa^{2}}
+\frac{\mathcal{F}(\mathcal{R})}{2\kappa^{2}},
\end{align}	
Using the value of $\sigma$, we get from (\ref{g7}) that
\begin{align}\label{g10}
p=\frac{(3\gamma-\mu)\mathcal{F}_\mathcal{R}(\mathcal{R})}{\kappa^{2}}-\frac{\mathcal{F}(\mathcal{R})}{2\kappa^{2}},
\end{align}	
This completes the proof.\par
Using Theorem \ref{th8} in Lemma \ref{l1}, we get the subsequent result:
\begin{corollary}\label{cor4.2}
For $\mathcal{R}=$ constant, in any $(QC)_4$ spacetime solution of the $\mathcal{F}(\mathcal{R})$-gravity, the effective pressure and density are described by
\begin{equation*}
  \sigma^{(eff)}=(3\mu -3\gamma+\mathcal{R})\frac{\mathcal{F}_\mathcal{R}(\mathcal{R})}{2\kappa^{2}},
\end{equation*}
\begin{equation*}
  p^{(eff)}=(6\gamma-2\mu-\mathcal{R})\frac{\mathcal{F}_\mathcal{R}(\mathcal{R})}{2\kappa^{2}}.
\end{equation*}
\end{corollary}
\section{Energy Conditions}
Now, we are dealing with a PF matter distribution and in accordance with the research by Capozziello et al. \cite{cap2}, the ECs retrieved from standard GR are
\begin{align}\label{5.5}
	NEC& \;\;\;if \;\;\;\sigma^{(eff)}+p^{(eff)}\geq0,\\\label{5.6}
	WEC& \;\;\;if \;\;\;\sigma^{(eff)}\geq0\quad\mathrm{and}\quad\sigma^{(eff)}+p^{(eff)}\geq0,\\\label{5.7}
	DEC& \;\;\;if \;\;\;\sigma^{(eff)}\geq0\quad\mathrm{and}\quad\sigma^{(eff)}\pm p^{(eff)}\geq0,\\\label{5.8}
	SEC& \;\;\;if \;\;\;\sigma^{(eff)}+3p^{(eff)}\geq0.
\end{align}
We now examine the ECs for a $\mathcal{F}\left(\mathcal{R}\right)$-gravity model in the following subsection.\par
{\bf A. Model:} $\mathcal{F}\left(\mathcal{R}\right)=e^{\mathcal{R}}\log\mathcal{R}
-\log\mathcal{R}-\mathcal{R}-\sum\limits_{l=2}^{\infty}\dfrac{\mathcal{R}^{l}}{l\cdot l!}$\par
For this model, with the help of equation \eqref{b2} and corollary $\ref{cor4.2}$, the effective energy density and pressure are presented as
\begin{align}\label{4.26}
\sigma^{(eff)}=\dfrac{1}{2\kappa^{2}}\left(9\gamma-3\mu\right)e^{\mathcal{R}}\log\mathcal{R},
\end{align}
\begin{equation}\label{4.27}
p^{(eff)}=\dfrac{1}{2\kappa^{2}}\left(4\mu-6\gamma\right)e^{\mathcal{R}}\log\mathcal{R}.
\end{equation}
The ECs for this setup can now be discussed using \eqref{4.26} and \eqref{4.27}.\par
\begin{tabulary}{\linewidth}{CC}
	
	\includegraphics[height=0.45\textheight]{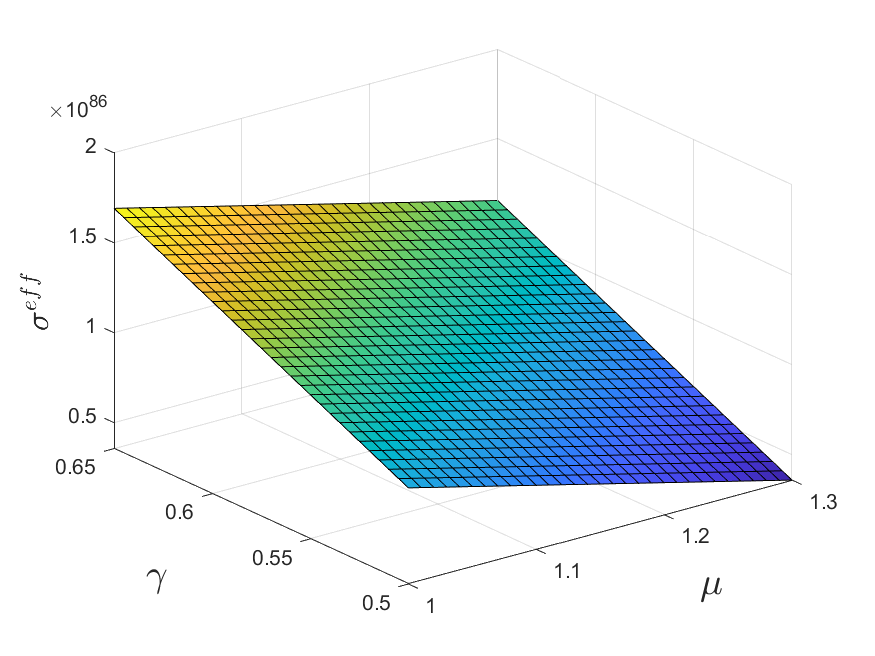}

	{\bf Fig. 1:} Evolution of energy density with respect to $\mu$ and $\gamma$
	
\end{tabulary}
\begin{tabulary}{\linewidth}{CC}
	
	\includegraphics[height=0.45\textheight]{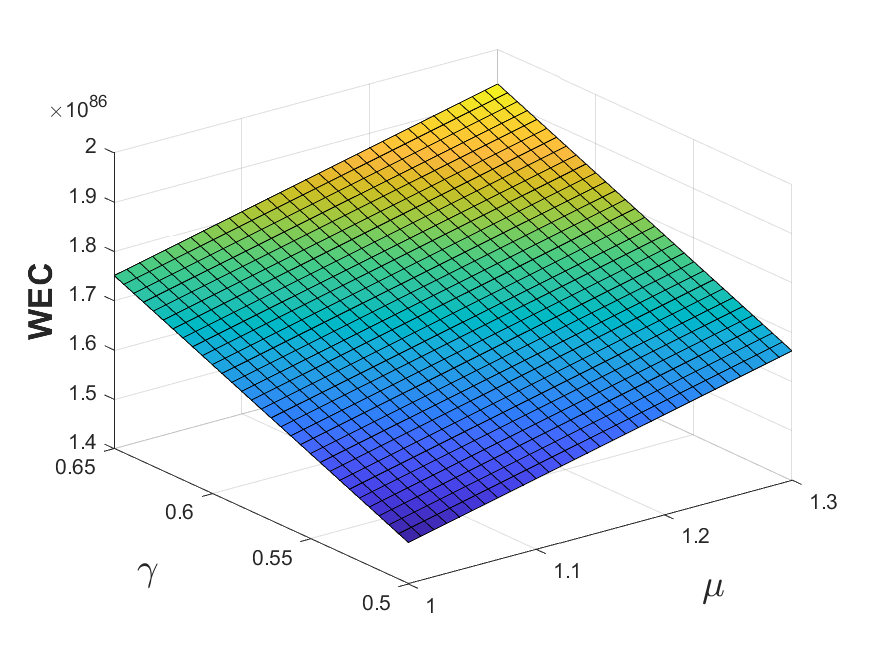}

	{\bf Fig. 2:} Evolution of WEC with respect to $\mu$ and $\gamma$
	
\end{tabulary}
Figs. $1$ and $2$, respectively, show the profiles of $\sigma^{eff}$ and WEC. The effective energy density is positive throughout the evolution of Ricci scalar $\mathcal{R}>1$, $\mu>0$ and $\gamma>\dfrac{\mu}{3}$\,. One can see from Fig. $1$ that for higher values of $\mu$, the effective energy density is high. Fig. 2 shows the WEC profile, which has a positive range for its value. As NEC is a part of WEC. Consequently, NEC and WEC are satisfied.
\begin{tabulary}{\linewidth}{CC}
	
	\includegraphics[height=0.45\textheight]{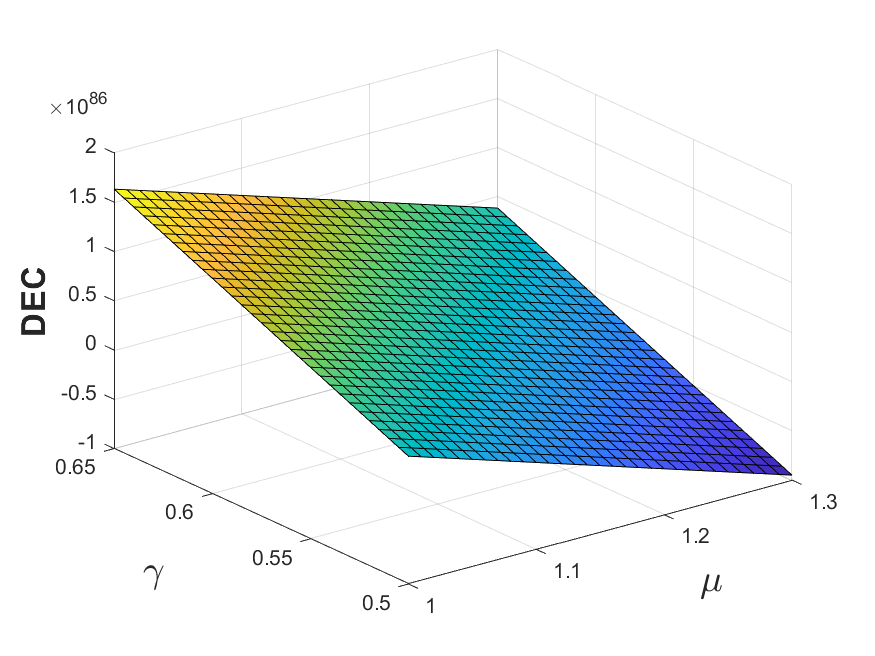}

	{\bf Fig. 3:} Evolution of DEC with respect to $\mu$ and $\gamma$
	
\end{tabulary}
\begin{tabulary}{\linewidth}{CC}
	
	\includegraphics[height=0.45\textheight]{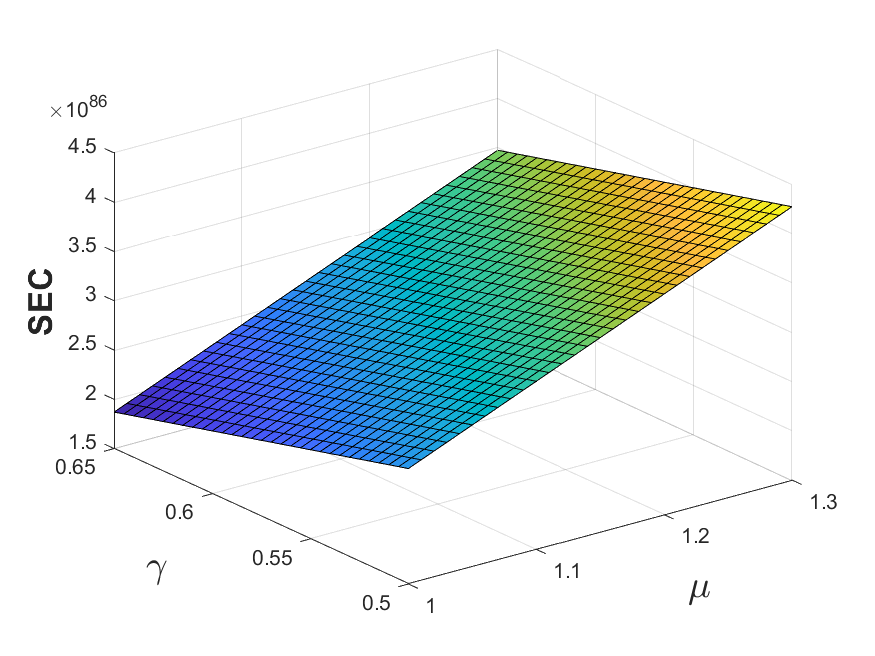}

	{\bf Fig. 4:} Evolution of SEC with respect to $\mu$ and $\gamma$
	
\end{tabulary}
To finish up our discussions, we have drawn DEC and SEC in the Figs. $3$ and $4$ for the above model. From these figures it is obvious that DEC does not hold but SEC does.
\section{Discussion}
Spacetime, a torsion-free, time-oriented Lorentzian manifold, is the stage on which the physical world is now being modelled. According to GR theory, the universe's matter content may be determined by choosing the appropriate energy momentum tensor, and is accepted to act like a PF- spacetime in the cosmological models.\par
Here, we study a spacetime of quasi- constant sectional curvature, and we have demonstrated that this spacetime is consistent with the universe's current condition and represents a PF- spacetime. Additionally, if the spacetime is Ricci symmetric or Ricci semisymmetric, it either represents a spacetime of constant sectional curvature or phantom era. We further establish that a Ricci symmetric $(QC)_{4}$-spacetime implies a static spacetime and is of Petrov type I, D, or O.\par

The investigation of $(QC)_{4}$-spacetime within the context of $\mathcal{F}\left(\mathcal{\mathcal{R}}\right)$-gravity has been the foremost concern of this article. Here, both analytic and graphical analysis of our investigations have been done. For a better understanding, we have applied the analytical method to develop our formulation and graphical analysis has been done to assess the stability of one cosmological toy model, such as $\mathcal{F}\left(\mathcal{R}\right)=e^{\mathcal{R}}\log\mathcal{R}
-\log\mathcal{R}-\mathcal{R}-\sum\limits_{l=2}^{\infty}\dfrac{\mathcal{R}^{l}}{l\cdot l!}$.\par
Additionally, we looked at the cosmological models' stability analysis using ECs. Figures 1, 2, 3, and 4 show the profiles of ECs for our model. The evolution of the energy density for parameters $\mu > 1$ and $\gamma > 0.5$ has been found to be positive. However, DEC broke the agreement, whereas SEC, NEC and WEC have been satisfied. The accelerated expansion of the cosmos is consistent with each of the aforementioned energy condition profiles. This agrees with the findings from the prior investigation of ECs in $\mathcal{F}\left(\mathcal{\mathcal{R}}\right)$-gravity.

\section{Declarations}
\subsection{Funding }
Not applicable.
\subsection{Conflicts of interest/Competing interests}
The authors declare that they have no conflict of interest and all authors contributed equally to this work.
\subsection{Availability of data and material }
Not applicable.
\subsection{Code availability}
Not applicable.

\end{document}